\documentclass[fleqn]{mat01}
\usepackage{times,mathtimy,amssymb}
\begin{document}

\setcounter{page}{1}
\firstpage{1}

\font\xx=msam5 at 10pt
\def\ab{\mbox{\xx{\char'03}}}

\font\bi=tibi at 7.3pt

\newtheorem{theo}{\bf Theorem}
\renewcommand\thetheo{\arabic{theo}}
\newtheorem{propo}[theo]{\rm PROPOSITION}
\newtheorem{coro}[theo]{\rm COROLLARY}
\newtheorem{lem}[theo]{Lemma}
\newtheorem{rema}[theo]{Remark}

\newcommand{\R}{\mbox{$\mathbb{R}$}}
\newcommand{\N}{\mbox{$\mathbb{N}$}}
\newcommand{\Q}{\mbox{$\mathbb{Q}$}}
\newcommand{\C}{\mbox{$\mathbb{C}$}}
\newcommand{\Z}{\mbox{$\mathbb{Z}$}}
\newcommand{\K}{\mbox{$\mathbb{K}$}}

\title{On some congruence with application to exponential sums}

\markboth{Soon-Mo Jung}{Study of an estimation of exponential sums}

\author{SOON-MO JUNG}

\address{Mathematics Section, College of Science and Technology, Hong-Ik
University, 339-701 Chochiwon, Korea\\
\noindent E-mail: smjung@wow.hongik.ac.kr}

\volume{114}

\mon{February}

\parts{1}

\Date{MS received 19 July 2003; revised 17 November 2003}

\begin{abstract}
We will study the solution of a congruence, $x \equiv
g^{(1/2)\omega_g(2^n)} \bmod 2^n$, depending on the integers $g$ and
$n$, where $\omega_g(2^n)$ denotes the order of $g$ modulo $2^n$.
Moreover, we introduce an application of the above result to the study
of an estimation of exponential sums.
\end{abstract}

\keyword{Congruence; multiplicative order; exponential sum.}

\maketitle

\section{Introduction}

Divisibility is an essential property between elements of algebraic
structures. As we know, the concept of congruence originates from that
of divisibility and it is not only a convenient notation but also a very
helpful method to prove many theorems in number theory. For this reason,
almost every textbook for number theory treats this topic extensively
(cf.~\cite{2,3}).

There are very close relationships between the congruences and the
exponential sums. In general, we may frequently use the properties of
congruences to estimate any exponential sum (see \cite{1,4,5}).

Discrepancy is also an important concept (or quantity) of number theory.
It measures the deviation of the sequence from an ideal distribution and
can be applied to problems in numerical analysis. We can easily
calculate, e.g., by the inequality of Erd\"{o}s--Tur\'{a}n, the upper
bound of the discrepancy if the related exponential sum has been
estimated. So, the estimations of the exponential sums are very
interesting.

Let $g$ and $m > 0$ be relatively prime integers. The (multiplicative)
{\em order\/} of $g$ modulo $m$ is defined as the least positive
exponent $k$ such that $g^k \equiv 1 \bmod m$, and we will denote it by
$\omega_g (m)$, i.e.,
\begin{equation*}
\omega_g (m) ~=~ \min\{ k \in \N : g^k \equiv 1 \bmod m \}.
\end{equation*}
In this paper, we investigate the solutions of a congruence related to
the order of an odd integer modulo $2^n$. These results will be applied
to the study of an estimation of exponential sums.

\section{Congruence modulo ${\bf 2}^{\hbox{{\bi n}}}$}

For a positive integer $n$ and an odd integer $g$, we recall that the
order of $g$ modulo $2^n$ is defined by
\begin{equation*}
\omega_g (2^n) ~=~ \min\{ k \in \N : g^k \equiv 1 \bmod 2^n \}.
\end{equation*}

\begin{lem} Given an $n \in \N${\rm ,} assume that $g$ is an odd
integer such that $g \not\equiv \pm 1 \bmod 2^n$. Then
\begin{equation*}
\omega_g (2^{n+1}) ~=~ 2\omega_g (2^n).
\end{equation*}
\end{lem}

\begin{proof} By the definition, we have $g^{\omega_g (2^n)} \equiv 1
\bmod 2^n$, i.e., $g^{\omega_g (2^n)} = a 2^n + 1$ for some $a \in \Z$.
Thus, for $k = 1,\ldots,\omega_g (2^n)-1$, we obtain
\begin{equation*}
g^k g^{\omega_g (2^n)} ~=~ a g^k 2^n + g^k ~\not\equiv~ 1 \bmod
2^{n+1},
\end{equation*}
since $g^k$ is odd and $g^k \not\equiv 1 \bmod 2^n$. (If $ag^k 2^n + g^k
\equiv 1 \bmod 2^{n+1}$, then there is an integer $b$ with $ag^k 2^n +
g^k = b 2^{n+1} + 1$ and hence $(ag^k - 2b) 2^n + g^k = 1$, i.e., $g^k
\equiv 1 \bmod 2^n$, a contradiction.)

However, it holds that
\begin{align*}
g^{\omega_g (2^n)} g^{\omega_g (2^n)} &= a g^{\omega_g (2^n)} 2^n +
g^{\omega_g (2^n)}\\
&= a (a 2^n + 1) 2^n + a 2^n + 1\\
&= a^2 2^{2n} + a 2^{n+1} + 1\\
&\equiv 1 \bmod 2^{n+1}.
\end{align*}
Therefore, it holds that $\omega_g (2^{n+1}) = 2\omega_g (2^n)$.
\hfill \ab
\end{proof}

\setcounter{theo}{0}
\begin{rema}
{\rm If $n \geq 2$ is an integer and $g$ is an odd integer with $g
\not\equiv \pm 1 \bmod 2^n$, then the previous lemma implies that
$\omega_g (2^n) = 2\omega_g (2^{n-1})$. If $g \equiv -1 \bmod 2^n$ then
$\omega_g (2^n) = 2$. Thus, for $n \geq 2$ and $g \not\equiv 1 \bmod
2^n$, $(1/2)\omega_g (2^n)$ is a positive integer.}
\end{rema}

The following lemma reveals that the congruence $x \equiv
g^{(1/2)\omega_g (2^n)} \bmod 2^n$ has at most three distinct solutions
(modulo $2^n$). However, as we see in Lemmas 3 and 4 below, the
congruence can have only two distinct solutions, $-1$ and $2^{n-1} + 1$
modulo $2^n$, indeed.

\setcounter{theo}{1}
\begin{lem}
Assume that $n$ is an integer $\geq 3$ and $g$ is an odd integer with $g
\not\equiv 1 \bmod 2^n$. Every solution of the following congruence
\begin{equation*}
x ~\equiv~ g^{(1/2)\omega_g (2^n)} \bmod 2^n
\end{equation*}
is congruent to one of $\{ -1, 2^{n-1} - 1, 2^{n-1} + 1 \}$ modulo
$2^n$.
\end{lem}

\begin{proof}
Since $g$ is odd, so is $g^{(1/2)\omega_{g}(2^{n})}$. Set
\begin{equation}
g^{(1/2)\omega_g (2^n)} ~=~ a 2^{n-1} + 2k + 1,\label{eq:sm1}
\end{equation}
where $a$ and $k$ are some integers with $0 \leq k < 2^{n-2}$. It then
follows from the definition of $\omega_g (2^n)$ that
\begin{align*}
g^{\omega_g (2^n)} &= (a 2^{n-1} + 2k + 1)^2\\
 &= 2^n (a^2 2^{n-2} + a (2k + 1)) + 4k^2 + 4k + 1\\
 &\equiv 1 \bmod 2^n .
\end{align*}
Therefore, we get $2^n|(4k^2 + 4k)$, i.e., $2^{n-2}|k(k+1)$.

\begin{enumerate}
\renewcommand\labelenumi{(\alph{enumi})}
\item Suppose $k$ is odd. Then, $k+1$ is even. Hence, it follows that
$2^{n-2}|(k+1)$. As $1 \leq k+1 \leq 2^{n-2}$, we conclude that
$k+1 = 2^{n-2}$. Due to (\ref{eq:sm1}), we have
\begin{equation*}
\hskip -1.25pc g^{(1/2) \omega_g (2^n)} + 1 ~=~ (a+1) 2^{n-1}.
\end{equation*}
If we square both sides of the last equality, then
\begin{equation*}
\hskip -1.25pc g^{\omega_g (2^n)} + 2g^{(1/2) \omega_g (2^n)} + 1 ~=~ (a+1)^2 2^{2n-2}.
\end{equation*}
As $g^{\omega_g (2^n)} \equiv 1 \bmod 2^n$, it follows from the last
equality that $2 + 2g^{(1/2)\omega_g (2^n)} \equiv 0 \bmod 2^n$, i.e.,
$1 + g^{(1/2)\omega_g (2^n)} \equiv 0 \bmod 2^{n-1}$. By (\ref{eq:sm1})
and the last congruence, we have $k+1 \equiv 0 \bmod 2^{n-2}$, i.e.,
there exists an integer $b$ such that $k = b 2^{n-2} - 1$. If $a+b$ is
even, it then follows from (\ref{eq:sm1}) that $g^{(1/2)\omega_g (2^n)}
\equiv -1 \bmod 2^n$. If $a+b$ is odd, we get $g^{(1/2)\omega_g (2^n)}
\equiv 2^{n-1} - 1 \bmod 2^n$.

\item Now, suppose $k$ is even. Then, it holds that $2^{n-2}|k$. As
$0 \leq k < 2^{n-2}$, it must be $k = 0$. Due to (\ref{eq:sm1}), we
obtain
\begin{equation*}
\hskip -1.25pc g^{(1/2)\omega_g (2^n)} ~=~ a 2^{n-1} + 1.
\end{equation*}
If $a$ is even, then $g^{(1/2)\omega_g (2^n)} \equiv 1 \bmod 2^n$ and
$(1/2)\omega_g (2^n) < \omega_g (2^n)$, a contradiction. Hence, $a$ is
odd and $g^{(1/2)\omega_g (2^n)} \equiv 2^{n-1} + 1 \bmod 2^n$.
\hfill \ab
\end{enumerate}\vspace{-1.3pc}
\end{proof}

In the previous lemma, we have proved that if $n \geq 3$ and $g$
($\not\equiv 1 \bmod 2^n$) is an odd integer, it may occur that
$g^{(1/2)\omega_g (2^n)} \equiv -1 \bmod 2^n$. In the following lemma,
we will prove that if $g \not\equiv -1 \bmod 2^n$ then $g^{(1/2)\omega_g
(2^n)} \not\equiv -1 \bmod 2^n$.

\begin{lem}
Assume that $n$ is an integer $\geq 3$ and $g$ is an odd integer with $g
\not\equiv 1 \bmod 2^n$. If the congruence
\begin{equation*}
g^{(1/2)\omega_g (2^n)} ~\equiv~ -1 \bmod 2^n
\end{equation*}
holds{\rm ,} then $\omega_g (2^n) = 2$ and $g \equiv -1 \bmod 2^n$.
\end{lem}

\begin{proof} Set $g = 2k+1$ with some $k \in \Z\!\setminus\!\{ 0 \}$.
By the binomial theorem, we have
\begin{align}
g^{(1/2)\omega_g (2^n)} &= (2k+1)^{(1/2)\omega_g (2^n)}\nonumber\\
&= \sum_{i=0}^{(1/2)\omega_g (2^n)}
\left( \begin{array}{c}
(1/2)\omega_g (2^n)\\[.2pc]
i
\end{array} \right)
(2k)^{i}.\label{eq:sra}
\end{align}
Assume that $g \not\equiv -1 \bmod 2^n$. It follows from Lemma~1 that
$\omega_g (2^n) = 2^2 \omega_g (2^{n-2})$. Hence,
\begin{equation*}
\frac{1}{2} \omega_g (2^n) ~=~ \left( \begin{array}{c}
(1/2)\omega_g (2^n)\\[.2pc]
1
\end{array} \right)
\end{equation*}
is even. Therefore, it follows from (\ref{eq:sra}) and the given
congruence that
\begin{align*}
g^{(1/2)\omega_g (2^n)} &= 1 + \sum_{i=1}^{(1/2)\omega_g (2^n)}
\left( \begin{array}{c}
(1/2)\omega_g (2^n)\\[.2pc]
 i
\end{array} \right) (2k)^i\\[.2pc]
&= 1 + 2^2 m \\
&\equiv -1 \bmod 2^n,
\end{align*}
where $m$ is an appropriate integer. Hence, $2^2 m \equiv -2 \bmod 2^n$
or $2^2 m = a 2^n - 2$ for some $a \in \Z$. Thus, we have $m = a 2^{n-2}
- 1/2$, which is a contradiction. Hence, it follows that $g \equiv -1
\bmod 2^n$. Moreover, we see that $\omega_g (2^n) = 2$.\hfill \ab
\end{proof}

In view of Lemma~2, we may expect that the congruence $g^{(1/2)\omega_g
(2^n)} \equiv 2^{n-1} - 1 \bmod 2^n$ holds for some $n \geq 3$ and some
odd integer $g$. However, we have to give up our expectation as we see
in the following lemma.

\begin{lem}
Assume that $n$ is an integer $\geq 3$ and $g$ is an odd integer with $g
\not\equiv 1 \bmod 2^n$. The congruence
\begin{equation*}
g^{(1/2)\omega_g (2^n)} ~\equiv~ 2^{n-1} - 1 \bmod 2^n
\end{equation*}
does not hold.
\end{lem}

\begin{proof} If we set $g = 2k+1$ with some $k \in \Z\!\setminus\!\{ 0
\}$ and follow the first part of the proof of Lemma 3, we then obtain
(\ref{eq:sra}).

\begin{enumerate}
\renewcommand\labelenumi{(\alph{enumi})}
\item Assume that $g \not\equiv -1 \bmod 2^n$. By following the second
part of the proof of Lemma~3, we have
\begin{equation*}
\hskip -1.25pc g^{(1/2)\omega_g (2^n)} ~=~ 1 + 2^2 m ~\equiv~ 2^{n-1} -1 \bmod 2^n,
\end{equation*}
where $m$ is some integer. Thus, $2^2 m \equiv 2^{n-1} - 2 \bmod 2^n$ or
$2^2 m = a 2^n + 2^{n-1} - 2$ for some $a \in \Z$. Hence, we get $m = a
2^{n-2} + 2^{n-3} - 1/2$, which is a contradiction.

\item If $g \equiv -1 \bmod 2^n$, then $\omega_g (2^n) = 2$. Therefore,
the given congruence reduces to $g \equiv 2^{n-1} - 1 \bmod 2^n$.
Altogether, it follows that $2^{n-1} - 1 \equiv -1 \bmod 2^n$, which is
a contradiction.
\hfill \ab
\end{enumerate}
\end{proof}\vspace{-1.4pc}

In the following theorem, we will summarize all the results of lemmas~2,
3 and 4.

\begin{theo}[\!]
Let $n$ be a given integer larger than $2$ and let $g$ be a given odd
integer.

\begin{enumerate}
\renewcommand\labelenumi{\rm (\alph{enumi})}
\item If $g \not\equiv \pm 1 \bmod 2^n${\rm ,} then
\begin{equation*}
\hskip -1.25pc g^{(1/2)\omega_g(2^n)} ~\equiv~ 2^{n-1} + 1 \bmod 2^n.
\end{equation*}

\item If $g \equiv -1 \bmod 2^n${\rm ,} then
\begin{equation*}
\hskip -1.25pc g^{(1/2)\omega_g(2^n)} ~\equiv~ -1 \bmod 2^n.
\end{equation*}
\end{enumerate}
\end{theo}

\section{Application to exponential sums}

For any non-zero integer $w$, we denote by $d(w)$ the {\em greatest
exponent\/} $k$ of $2$ satisfying $2^k|w$, i.e., $d(w) = \max\{ k
\in \N_0 : 2^k|w \}$, and we use the notation $2^{d(w)}\|
w$ for this\break case.

\setcounter{theo}{1}
\begin{rema}
{\rm Every non-zero integer $w$ can be represented by $w = 2^{d(w)} w_0$,
where $w_0$ is an odd integer.}
\end{rema}

For an odd integer $g \not\in \{-1,1\}$, we define
\begin{equation*}
c(g) ~=~ \begin{cases}
 \min\{ k \in \N : g < 2^{k-1} - 1 \} &\mbox{for}~ g > 1,\\
 \min\{ k \in \N : g > -2^{k-1} - 1 \} &\mbox{for}~ g < -1.
 \end{cases}
\end{equation*}
In this section, we will introduce an application of the previous
results. Indeed, Korobov \cite{korobov} has already proved a similar
result in the following theorem. However, our proof is more visible than
that of Korobov.

\setcounter{theo}{5}
\begin{theo}[\!]
If $g \not\in \{ -1, 1 \}$ is an odd integer and $w$ is a non-zero
integer{\rm ,} then
\begin{equation*}
\sum_{k=1}^{\omega_g (2^n)} {\rm e}^{2 \pi iwg^k / 2^n} = 0
\end{equation*}
for any integer $n \geq d(w) + \max\{ 3, c(g) \}$.
\end{theo}

\begin{proof}$\left.\right.$

\begin{enumerate}
\renewcommand\labelenumi{\rm (\alph{enumi})}
\item Assume that $w$ is odd. As $n \geq c(g)$, if $g > 1$ then $1 < g <
2^{n-1} - 1$. Similarly, if $g < -1$ then $2^{n-1} - 1 < 2^n + g < 2^n -
1$. Hence, we conclude that $g \not\equiv \pm 1 \bmod 2^n$.

It follows from Theorem~5(a) that
\begin{equation*}
\hskip -1.25pc g^{(1/2)\omega_g (2^n)} ~\equiv~ 2^{n-1} + 1 \bmod 2^n
\end{equation*}
or
\begin{equation}
\hskip -1.25pc g^{(1/2)\omega_g (2^n)} ~=~ a 2^n + 2^{n-1} + 1\label{eq:sm2}
\end{equation}
for some $a \in \Z$.

Since $g^k$ is odd, (\ref{eq:sm2}) yields
\begin{align*}
\hskip -1.25pc g^{k + (1/2)\omega_g (2^n)} &= a g^k 2^n + g^k 2^{n-1} + g^k \\
\hskip -1.25pc &\equiv 2^{n-1} + g^k \bmod 2^n
\end{align*}
for $k = 1, \ldots, (1/2)\omega_g (2^n)$. Moreover, since $w$ is assumed
to be odd, we get
\begin{align}
\hskip -1.25pc w g^{k + (1/2)\omega_g (2^n)} &\equiv w 2^{n-1} + w g^k \bmod
2^n\nonumber\\
\hskip -1.25pc &\equiv 2^{n-1} + w g^k \bmod 2^n.\label{eq:sm3}
\end{align}
Hence, it follows from (\ref{eq:sm3}) that
\begin{align*}
\hskip -1.25pc \sum_{k=1}^{\omega_g (2^n)} {\rm e}^{2\pi iwg^k / 2^n} &=
\sum_{k=1}^{(1/2)\omega_g (2^n)} ({\rm e}^{2\pi iwg^k / 2^n} + {\rm
e}^{2\pi iwg^{k + (1/2)\omega_g (2^n)} / 2^n})\\[.2pc]
\hskip -1.25pc &= \sum_{k=1}^{(1/2)\omega_g (2^n)} ({\rm e}^{2\pi iwg^k / 2^n} + {\rm
e}^{2\pi i(2^{n-1} + w g^k) / 2^n})\\
\hskip -1.25pc &= 0.
\end{align*}

\item Assume now that $w$ is a non-zero even integer. As we noticed in
Remark~2, $w$ may be represented as $w = 2^{d(w)} w_0$, where $w_0$ is
odd. If we set $m = n - d(w)$, we obtain $g \not\equiv \pm 1 \bmod 2^m$
by following the first part of (a).

It follows from Theorem~5(a) that
\begin{equation*}
\hskip -1.25pc g^{(1/2)\omega_g (2^m)} ~\equiv~ 2^{m-1} + 1 \bmod 2^m .
\end{equation*}
By a similar method given in (a), we get, instead of (\ref{eq:sm3}),
\begin{equation*}
\hskip -1.25pc w_0 g^{k + (1/2)\omega_g (2^m)} ~\equiv~ 2^{m-1} + w_0 g^k \bmod 2^m
\end{equation*}
for $k = 1, \ldots, (1/2)\omega_g (2^m)$. As we did in the last part of
(a), we have
\begin{equation}
\hskip -1.25pc \sum_{k=1}^{\omega_g (2^m)} {\rm e}^{2\pi iw_0 g^k / 2^m} ~=~
0.\label{eq:sm4}
\end{equation}
According to Lemma~1, we get
\begin{equation}
\hskip -1.25pc \omega_g (2^n) ~=~ 2^{d(w)} \omega_g (2^m).\label{eq:sm5}
\end{equation}
Finally, it follows from (\ref{eq:sm4}) and (\ref{eq:sm5}) that
\begin{align*}
\hskip -1.25pc \sum_{k=1}^{\omega_g (2^n)} {\rm e}^{2\pi iwg^k / 2^n} &=
\sum_{k=1}^{\omega_g (2^n)} {\rm e}^{2\pi iw_0 g^k / 2^m}\\
\hskip -1.25pc &= 2^{d(w)} \sum_{k=1}^{\omega_g (2^m)} {\rm e}^{2\pi iw_0 g^k / 2^m}\\
\hskip -1.25pc &= 0,
\end{align*}
as required.\hfill \ab
\end{enumerate}
\end{proof}\vspace{-1.5pc}

\section*{Acknowledgement}

This work was supported by the 2004 Hong-Ik University Academic Research
Support Fund.

\end{document}